\documentclass[final,12pt,english]{article}
\usepackage{verbatim}

\usepackage{latexsym}
\usepackage{amsmath,amsthm,amsfonts}
\usepackage{eucal}
\usepackage[notref,notcite]{showkeys}

\usepackage{mathtext}
\usepackage[cp1251]{inputenc}
\usepackage[T2A]{fontenc}

\usepackage[dvips]{graphicx}
\usepackage{amsmath}
\usepackage{amssymb}
\usepackage{amsxtra}
\usepackage{amsthm}
\usepackage{ifthen}

\newtheorem{theorem}{Theorem}[section]
\newtheorem{definition}{Definition}[section]
\newtheorem{lemma}{Lemma}[section]

\newtheorem{remark}{Remark}[section]

\numberwithin{equation}{section}

\DeclareMathOperator{\supp}{supp}

\begin{document}

\title {Euler equation for incompressible non-Newtonian fluids: finite speed of propagations
and asymptotic behavior of weak solutions }

\author{{\normalsize\bf Yuliya V. Namlyeyeva and Roman M. Taranets}
\footnote{Research is partially supported by the INTAS project
Ref. No: 05-1000008-7921}
\smallskip}

\maketitle

\begin{abstract}
We investigate multidimensional model for incompressible
\linebreak non-Newtonian fluids. Using method of energy estimates
we prove the property of finite speed of propagations of the
solution support for this problem. We find sharp bounds of the
propagations by $L^2$--norm and $L^1$--norm of initial data.
\end{abstract}


\textbf{keywords:} non-Newtonian fluids, asymptotic behavior,
finite speed of propagations of support

\section{Formulation of the problem and main results}
There are many substances like geological materials, liquid foams,
polymeric fluids, which are capable of flowing but which exhibit
flow characteristics that cannot be adequately described by the
classical linearly viscous fluid model. In order to describe some
of the departures from Newtonian behavior evinced by such
materials, many idealized material models have been suggested.
There are several phenomena which appear when studying
non-Newtonian fluids: the shear thinning and the shear thickening,
the ability of a creep, the ability to relax stresses, the
presence of normal stress differences in simple shear flow, the
presence of yield stress. For more details see \cite {Ra}. The
equations governing the unsteady motion these fluids have been
previously investigated, e.g. \cite{BBN,MNRR,Po1,Po2}, etc.

Let components of the symmetric deformation velocity tensor are
given by
\begin{equation}
({\bf D}u)_{ij} = \frac12 \Big( \frac{\partial u_i}{\partial x_j}
+ \frac{\partial u_j}{\partial x_i}\Big)
\end{equation}
$(|{\bf D} u|=({\bf D}u_{ij} {\bf D}u_{ij})^{\frac12})$.

We consider  multi--dimensional problem  of non-Newtonian model
which has the following form
\begin{gather}\ u_t + (u\cdot \nabla )\,u =
 \mu_1 div(|\mathbf{D} u|^{p - 2} \mathbf{D}u)
+ \nabla \pi, \hfill \label{0} \\
div \, u= 0, \hfill \label{1}\\
u(0,x) = u_0(x).\  \hfill \label{2}
\end{gather}
\vspace{-2.5cm}
$$
(C)\qquad\left\{ \hspace{+3cm}\phantom {\begin{gathered} \ u_t +
(u\cdot \nabla )\,u = \mu_0 \triangle u + \mu_1 div (|D u|^{p -
2}D u), \hfill \nonumber \\
div \, u= 0, \hfill \label{1}\\
u(0,x) = u_0(x).\  \hfill \label{2}
\end{gathered}} \right.
$$
Here $u$ is the velocity field, $\pi$ is the pressure for an
incompressible power-law fluids, $u_0$ is the initial value of the
velocity. Let $\mu_1
> 0,$
\begin{align}\label {supp}
u_0\in (W^{1,2}(\mathbb{R}^N_{-}))^N \cap H,\ \supp \,
u_{0,k}(x)\subseteq \mathbb{R}^N_{-},\ k=\overline{1,N},
\end{align}
where $ \mathbb{R}^N_{-}:=\{x=(x',x_N)\in \mathbb{R}^N: x_N < 0\}$
and $H=\{u\in L_2(\mathbb{R}^N)^N: \ div \, u=0\}$.

Existence of a weak solution for $p>\frac{3N}{N+2}$ and the
uniqueness and regularity for $p\geqslant \frac{3N + 2}{N+2}$,
$N=2,3$ were proved in \cite{BBN,MNRR}. The more detailed
information about qualitative behavior of the solution can be
found in \cite{A}.

\medskip

\begin{definition}
The function $u(x,t)=(u_1(t,x),\ldots , u_N(t,x))$ such that
\begin{align}\nonumber
&u\in L_p(0,T;W^{1,p}(\mathbb{R}^N)^N) \cap
L_2(0,T;W^{1,2}(\mathbb{R}^N)^N)\cap C(0,T;H), \\
&u_t \in L_2(0,T;H)\nonumber
\end{align}
is called a weak solution to problem $(C)$ if for a.e. $t > 0$,
the integral identity
\begin{equation}\label{el}
\int\limits_{\mathbb{R}^N} {u_t (t) \varphi \, dx} +
\int\limits_{\mathbb{R}^N} (u\cdot \nabla )\,u \,\varphi \,dx +
\mu_1 \int\limits_{\mathbb{R}^N} { |\mathbf{D} u(t)|^{p - 2}
\mathbf{D}u(t)\,: \mathbf{D}\varphi \, dx} = 0
\end{equation}
is satisfied for every $\varphi  \in
L_p(0,T;W^{1,p}(\mathbb{R}^N)^N) \cap
L_2(0,T;W^{1,2}(\mathbb{R}^N)^N)$ with $\text{div}\, \varphi = 0$.
\end{definition}

\smallskip

Basic definitions and properties of $p$--Laplacian we can find in
the books of Vazquez \cite{Va1,Va2}. The evolution $p$--Laplacian
is one of the most widely researched equations in the class of
nonlinear degenerate parabolic equations, already studied by
Raviart \cite{Rav}. Large time behaviour of the solutions was
investigated in Kamin, Vazquez \cite{KV} and they proved that
explicit solutions found by Barenblatt in 1952 are essentially the
only positive solutions to the Cauchy problem with initial data
$$u(x,0)=M\delta (x),\ M>0.$$ They established that any
nonnegative solution with globally integrable initial values is
asymptotically equal to the Barenblatt solution as $t \to \infty
$. In work of Lee, Petrosyan and Vazquez \cite {LPV} the property
of asymptotic concavity was proved.

Now we give brief explanations of our method of proof of the
finite speed of propagations. This method we call as the method of
nonhomogeneous functional inequalities and it is some adaptation
the energy method. It is connected with nonhomogeneous variants of
Stampacchia lemma, in fact, it is an adaptation to higher order
equations of local energy or Saint--Venant principle like
estimates method. It was elaborated with respect to higher order
quasilinear parabolic equations of the monotone type in
\cite{SS,SSH,SH2}. Also it was developed with regard to thin-film
equations with nonlinear convection in \cite{T5,T3}. Fluid
mechanics is one of the most natural fields for the application of
energy methods (see, \cite{A,AO1,AO2}). This is because the
fundamental conservation and balance laws employed for its
description suggest an adequate choice of energy functions which,
in turn, produces the formation of a free boundary. "Energy
methods are special interest in those situations in which
traditional methods based a comparison principles have failed.
This method yields the formation of a free boundary; in other
words, this means that the support of solution is localized in
space-time domain."(see,\cite{A})

\medskip

The main results are the following.
\medskip

\begin{theorem}\label{Th1}
Let $u(x,t)$ be a weak solution of problem $(C)$. Let $p \geqslant
\tfrac{3N + 2}{N + 2}$. Then there exists the function $\Gamma(t)
\in C[0, T] ,\ \Gamma(0) = 0$ such that
\begin{equation}\label{est}
\Gamma(t) =  c_1\,\max \{ t^{\frac{2}{2p + N(p- 2)}},\ t^{\frac{2p
+ N(p- 3)}{2p + N(p- 2)}}\}\ \forall\, t >0
\end{equation}
and
\begin{equation}\label{fsp}
\supp\,u_k(t,.) \cap \{x=(x',x_N) \in \mathbb{R}^N: x_N \geqslant
\Gamma(t)\} = \emptyset, \, k=\overline{1,N},
\end{equation}
where $c_1= c_1(p, \mu_1, \|u_0\|_{L_2(\mathbb{R}^N)^N})$ is a
positive constant depending on known parameters only.
\end{theorem}

\medskip

\begin{remark}
If (\ref{fsp}) is valid then a solution to $(C)$ has the property
of finite speed of propagation of the solution support.
\end{remark}

\smallskip

\begin{remark}
The estimate (\ref{est}) is sharp for scalar equation if the
initial function $u_0(x)$ belongs $L_2(\mathbb{R}^N)^N$ (see, for
example, \cite{SS}).
\end{remark}

\medskip

\begin{theorem}\label{Th2}
Let $u(x,t)$ be a weak solution of problem $(C)$. Let
$\|u\|_{L_1(\mathbb{R}^N)^N} \leqslant
\|u_0\|_{L_1(\mathbb{R}^N)^N}$ and $p \geqslant \tfrac{3N + 1}{N +
1}$. Then the function $\Gamma(t)$ has the following form
\begin{equation}\label{eest}
\Gamma(t) =  c_2\,\max \{ t^{\frac{1}{p + N(p - 2)}},\,t^{\frac{p
+ N(p - 3)}{p + N(p - 2)}}\}\ \forall\, t >0,
\end{equation}
where $c_2 = c_2(p, \mu_1, \|u_0\|_{L_1(\mathbb{R}^N)^N})$ is
positive constant.
\end{theorem}

\smallskip

\begin{remark}
The exponent $1/(p + N(p - 2))$ from (\ref{eest})  is well-known
Barenblatt exponent for $p$--Laplacian scalar equation (see, for
example, \cite{AT}).
\end{remark}

\section{Proof of Theorem 1.1}
\begin{lemma}
Let $u(x,t)$ be an arbitrary weak solution of problem $(C)$. Let
 $p \geqslant \tfrac{3N + 2}{N +2}$. Then the following estimate is valid
\begin{multline}\label{enq3}
\mathop {\sup }\limits_{t \in (0,T)} \int\limits_{\Omega (s +
\delta )} {|u(t)|^{2}\,dx}  + \frac{1} {T}\iint\limits_{Q_T (s +
\delta )} {|u|^{2}\,dx\,dt}  + \mu_1 \iint\limits_{Q_T (s +
\delta )} {|\nabla u|^{p}\,dx\,dt} \leqslant\\
\frac{c}{\delta^p} \iint\limits_{Q_T (s)}{|u|^{p}\,dx\,dt}+ \frac
{c}{\delta} \iint\limits_{Q_T (s)}{|u|^{3}\,dx\,dt},
\end{multline}
$\forall\, s \geqslant s_0 \geqslant 0,\ \delta > 0,\ T > 0$.
\end{lemma}

\begin{proof}
For an arbitrary $s \in \mathbb{R}^1$ and $\delta > 0$ we consider
the families of
 sets
$$
\begin{gathered}
\Omega (s) = \{ x = (x', x_N)\in \mathbb{R}^N: \ x_N \geqslant s
\},\ Q_T(s)= (0,T) \times \Omega (s) , \hfill\\
K (s,\delta ) = \Omega (s)\backslash \Omega (s + \delta ),\ K_T
(s,\delta ) = (0,T) \times  K (s,\delta ). \hfill
\end{gathered}
$$
Next we introduce our main sequence cut-off functions $\eta
_{s,\delta } (x) \in C^1 (\mathbb{R}^N)$, which possess the
following properties:
\begin{equation}\label{con1}
0 \leqslant \eta_{s,\delta } (x) \leqslant 1 \  \forall  x \in
\mathbb{R}^N, \ \eta _{s,\delta } (x) = \left\{
\begin{aligned} \hfill  0 \;
& , x \in \mathbb{R}^N \setminus \Omega(s),\\
\hfill  1\; & , x \in \Omega(s + \delta), \\
\end{aligned}\right.
\end{equation}
\begin{equation}\label{con11}
 \ | \nabla\,\eta _{s,\delta }| \leqslant \frac{c_5} {\delta }\
\forall\,x \in K (s,\delta ).
\end{equation}
Test integral identity (\ref{el}) by
\begin{equation}\label{z}
\varphi(x,t) = u(t,x)\psi(t,x),
\end{equation}
where $u(t,x)$ is a solution of the problem ($C$) and, without
loss of generality, we suppose that $(\nabla_x \psi(t,x), u(t,x))
= 0$. Here
$$
\psi(t,x)=\eta _{s,\delta }(x) \exp \left( { - t \cdot T^{ - 1} }
\right)\ \forall\,T > 0,
$$
Then we obtain
\begin{multline}\label{enq0}
\frac{1}{2} \int\limits_{\mathbb{R}^N} {\frac{\partial }{\partial
t}(u^{2}(t)) \psi\,dx} + \mu_1 \int\limits_{\mathbb{R}^N} {
|\mathbf{D} u(t)|^{p} \psi \,dx}+\\
\mu_1\int\limits_{\mathbb{R}^N} { |\mathbf{D} u(t)|^{p-2}\,(
\mathbf{D} u(t)\,\nabla \psi) \, u(t) \, dx}+
 \int\limits_{\mathbb{R}^N} (u\cdot \nabla )u(t)\, u(t)\psi
dx= 0.
\end{multline}
We consider integral terms in the last identity. After using
Young's inequality we get
\begin{multline}\label{in1}
\ \int\limits_{\mathbb{R}^N} { |\mathbf{D} u|^{p-2}\,( \mathbf{D}
u\,\nabla \psi) \, u \, dx}\leqslant \varepsilon_1\exp \left( { -
t \cdot T^{ - 1} } \right)\ \int\limits_{K(s,\delta )}{|\mathbf{D}
u|^{p}\,dx} + \\
\frac{c(\varepsilon_1)}{\delta^p} \exp \left( { - t \cdot T^{ - 1}
} \right)\ \int\limits_{K (s,\delta )}{|u|^{p}\,dx},
\end{multline}
\begin{multline}\label{in3}
\int\limits_{\mathbb{R}^N} (u\cdot \nabla )\,u \, u(t)\psi dx =
\frac12 \exp \left( { - t \cdot T^{ - 1} } \right)\
\int\limits_{\mathbb{R}^N}|u(t)|^2 \, u(t)\, \nabla \psi dx
\leqslant\\
 \frac 1{2\delta}\exp \left( { - t \cdot T^{ - 1} }
\right)\ \int\limits_{\mathbb{R}^N}|u|^3dx,
\end{multline}
Integrating (\ref{enq0}) respect with time, and taking into
account (\ref{in1}) and (\ref{in3}), after simple computation, we
deduce that
\begin{multline}\label{enq1}
\mathop {\sup }\limits_{t \in (0,T)} \int\limits_{\Omega (s +
\delta )} {|u(t)|^{2}\,dx}  + \frac{1} {T}\iint\limits_{Q_T (s +
\delta )} {|u|^2\,dx\,dt} +2\mu_1 \iint\limits_{Q_T (s + \delta )}
{|\textbf{D} u|^{p}\,dx\,dt} \leqslant \\
 \int\limits_{\Omega (s)} {|u_0(x)|^{2}\,dx}+
\frac{c(\varepsilon_1)}{\delta^p} \iint\limits_{K_T (s,\delta
)}{|u|^{p}\,dx\, dt}+ \frac{c}{\delta}\iint\limits_{K_T (s,\delta
)}|u|^3 dx\,dt + \mu_1\varepsilon_1\iint\limits_{K_T (s, \delta )}
{|\textbf{D} u|^{p}\,dx\,dt},
\end{multline}
where $s \in \mathbb{R}^1 ,\ \delta
> 0,\ T > 0$.
It follows from (\ref{supp}) that
\begin{equation}\label{id}
\int\limits_{\Omega (s)} {|u_0(x)|^{2}\,dx} = 0 \ \forall\,s
\geqslant s_0 \geqslant 0.
\end{equation}
Finally, from (\ref{enq1}), (\ref{id}) we obtain
\begin{multline}\label{in5}
\mathop {\sup }\limits_{t \in (0,T)} \int\limits_{\Omega (s +
\delta )} {|u(t)|^{2}\,dx}  + \frac{1} {T}\iint\limits_{Q_T (s +
\delta )} {|u|^2\,dx\,dt} +2\mu_1 \iint\limits_{Q_T (s + \delta )}
{|\textbf{D} u|^{p}\,dx\,dt} \leqslant \\
 \frac{c(\varepsilon_1)}{\delta^p} \iint\limits_{K_T (s,\delta
)}{|u|^{p}\,dx\, dt}+ \frac{c}{\delta}\iint\limits_{K_T (s,\delta
)}|u|^3 dx\,dt + \mu_1\varepsilon_1\iint\limits_{K_T (s, \delta )}
{|\textbf{D} u|^{p}\,dx\,dt},
\end{multline}
$\forall\, s \geqslant s_0 \geqslant 0,\ \delta
> 0,\ T > 0$. Choosing $\varepsilon_1
> 0$ sufficiently small, and iterating the limit inequality (\ref{in5}), we get
(\ref{enq3}).
\end{proof}

\begin{proof}[Proof of Theorem \ref{Th1}]
We denote
$$R_T(s, \delta):=\frac{c}{\delta^p}
\iint\limits_{Q_T (s)}{|u|^{p}\,dx\,dt}+ \frac {c}{\delta}
\iint\limits_{Q_T (s)}{|u|^{3}\,dx\,dt}.$$
 We introduce the energy
functions related to our solution:
$$
A_T (s): = \iint\limits_{Q_T (s)} {|u|^{p}\,dx\,dt},\ B_T (s): =
\iint\limits_{Q_T  (s)} {|u|^3 \,dx\,dt }.
$$
Applying Nirenberg-Gagliardo's (see, Lemma A.2 of Appendix~A) and
H\"{o}lder's inequalities, we get
$$
\begin{gathered}
A_T (s + \delta) \leqslant c\,T^{\alpha_1} (R_T(s,\delta ))^{1 +
\beta_1},\ \alpha_1 = \tfrac{2p} {2p + N(p- 2)}, \ \beta_1 =
\tfrac{p(p - 2)} {2p + N(p- 2)},\hfill\\
B_T (s + \delta) \leqslant c\,T^{\alpha_2} (R_T(s,\delta ))^{1 +
\beta_2},\ \alpha_2 = \tfrac{2p + N(p - 3)}{2p + N(p- 2)}, \
\beta_2 = \tfrac{p}{2p + N(p- 2)}. \hfill\\
\end{gathered}
$$
Next we define the functions
$$
C_T(s) := (A_T (s))^{1 + \beta_2}+ (B_T (s))^{1 + \beta_1}.
$$
Then
\begin{equation}\label{D1}
C_T (s+\delta) \leqslant \tilde{c}\,F(T)\,(\delta^{-p\beta}C_T^{1
+ \beta_1}(s) + \delta^{-\beta} C_T^{1 + \beta_2}(s) ),
\end{equation}
where
$$
\beta = (1 + \beta_1)(1 + \beta_2),\ F(T)= \max \{ T^{\alpha_1(1 +
\beta_2)},\ T^{\alpha_2(1 + \beta_1)}\}.
$$
Now we choose the parameter $\delta > 0$ which was arbitrary up to
now. First, we introduce the notations
$$
\begin{gathered}
\delta _T^{(1)} (s): = \left[ {2\tilde{c}\,F(T)\,C_T^{\beta_1}
(s)} \right]^{\frac{1} {p\beta}},\ \delta _T^{(2)} (s): = \left[
{2\tilde{c}F(T)\,C_T^{\beta_2} (s)}\right]^{\frac{1}
{\beta}},\hfill\\
J_T (s): = \max \{ \delta _T^{(1)} (s),\delta _T^{(2)} (s)\} .
\hfill\\
\end{gathered}
$$
We obtain the following main functional relation for the functions
$J_T (s)$
\begin{equation}\label{J}
J_T (s + J_T (s)) \leqslant \varepsilon \, J_T (s)\ \forall \, s
\geqslant s_0 \geqslant 0, \ 0 < \varepsilon  < 1.
\end{equation}
We set $s =  - 2\delta ,\  \delta  = s' > 0$ in (\ref{enq1}) and
pass to the limit as $s' \to + \infty $. Using the boundedness of
functions $A_T (s)$ and $B_T (s)$, we get
\begin{equation}\label{J0}
C_T (0) \leqslant c\,F(T)\|u_0\|_{L_2}^{\beta}.
\end{equation}
Now we apply Lemma~A.1 to the function $J_T (s)$ of (\ref{J}). As
a result, we get
\begin{equation}\label{J1}
J_T (s) \equiv 0 \ \forall\,s \geqslant s_0+\tfrac{1}{1 -
\varepsilon} J_T (s_0).
\end{equation}
Let $s_0 = 0$. Then, in view of (\ref{J0}), we find
$$
J_T (0) \leqslant c \, \max \{F(T)^{\frac{1}{p(1 + \beta_2)}},\
F(T)^{\frac{1}{1 + \beta_1}} \}\leqslant c\,\max \{ T^{\frac{2}{2p
+ N(p- 2)}},\ T^{\frac{2p + N(p- 3)}{2p + N(p- 2)}}\}
$$
$\forall\,T > 0,\ c = c(p,\,\mu_1,\,N,\,\|u_0\|_{L_2})$. Choosing
in (\ref{J1})
$$
s = \Gamma (T) = c\,\max \{ T^{\frac{2}{2p + N(p- 2)}},\
T^{\frac{2p + N(p- 3)}{2p + N(p- 2)}}\},
$$
where $p \geqslant \tfrac{3N + 2}{N + 2}$, we obtain that
$J_T(\Gamma(T)) = 0$. Thus $u_k(T,x) \equiv 0$ for all $x \in
\{x=(x',x_N) : x_N \geqslant \Gamma(T)\},\ k=\overline{1,N}$. And
Theorem~1.1 is proved completely.
\end{proof}

\section{Proof of Theorem~1.2}

\begin{lemma}
Let $u(x,t)$ be an arbitrary weak solution of problem $(C)$. Let
$\|u\|_{L_1(\mathbb{R}^N)^N} \leqslant
\|u_0\|_{L_1(\mathbb{R}^N)^N}$ and $p \geqslant \tfrac{3N + 1}{N +
1}$. Then the following estimates for the decay rate are valid
\begin{equation}\label{esde}
A_T (s)+ B_T (s)\leqslant \tilde{c}({\Theta}) \, T \,\bigl( s^{ -
N(p - 1)} +s^{ -\frac{ 2N}{p + N (p - 3)}} \bigr)
\end{equation}
$\forall\, s > 0,\ T > 0$, where $\tilde{c}({\Theta}) =
\tilde{c}(p,N,{\Theta}),\ {\Theta} = \left\| {u_0 } \right\|_{L_1
(\mathbb{R}^N )^N}$.
\end{lemma}

\begin{proof} Applying the interpolation inequality of Lemma~A.2
in the domain $K(s,\delta )$ to the function $v = |u|$ for $a = d
= p,\ b = 1$, $d_2= c\,\delta^{- N(p - 1)/p}$, and integrating the
result with respect to time from $0$ to $T$, we obtain
\begin{multline}\label{es1}
\iint\limits_{K_T (s,\delta )} {|u|^{p} }dx\,dt \leqslant c\,
\delta ^{ - N(p - 1)}\, T\,{\Theta}^{p}
 + c\, T^{1 - \theta_1 } {\Theta}^{p(1 - \theta_1 )}\times \\
\times\Biggl( {\;\iint\limits_{K_T (s,\delta )} {|\nabla u|^p
dx\,dt }} \Biggr)^{\theta_1 }, \text{ where } \theta_1  =
\tfrac{N(p - 1)}{p + N(p - 1)}.
\end{multline}
Similarly, applying the interpolation inequality of Lemma~A.2 in
the domain $K(s,\delta )$ to the function $v = |u|$ for $a =3,\ d
= p,\ b = 1$, $d_2= c\,\delta^{- 2N/3}$, and integrating the
result with respect to time, we find that
\begin{multline}\label{es2}
\iint\limits_{K_T (s,\delta )} {|u|^{3}dx\,dt } \leqslant c\,
\delta ^{ - 2N}\, T\,{\Theta}^{3}
 + c\, T^{1 - \tfrac{3\theta_2 }{p}} {\Theta}^{3(1 - \theta_2 )}\times \\
\times\Biggl( {\;\iint\limits_{K_T (s,\delta )} {|\nabla u|^p
dx\,dt}} \Biggr)^{\tfrac{3}{p}\theta_2 }, \text{ where } \theta_2
= \tfrac{2\,N\,p }{3(p + N(p - 1))},\ p > \tfrac{3N}{N + 1}.
\end{multline}
Inserting (\ref{es1}) and (\ref{es2}) in (\ref{in5}) and applying
Young's $''\varepsilon''$--inequality, we get
\begin{multline}\label{enq22}
L_T(s + \delta):=\mathop {\sup }\limits_{t \in (0,T)}
\int\limits_{\Omega (s + \delta )} {u^{2}(t)\,dx}  + \frac{1}
{T}\iint\limits_{Q_T (s + \delta )} {u^{2}\,dx\,dt}  + \mu_1 C
\iint\limits_{Q_T (s + \delta )} {|\nabla u|^{p}\,dx\,dt}
\leqslant \\
\varepsilon_4\, \mu_1 \iint\limits_{K_T (s,\delta )}{|\nabla
u|^{p}\,dx\,dt}  + \varepsilon_3 \mathop {\sup }\limits_{t \in
(0,T)} \int\limits_{K(s, \delta)}{|u(t)|^{2}\,dx} +
\frac{\varepsilon_3}{T} \iint\limits_{K_T(s,
\delta)}{|u|^{2}\,dx\, dt} +\\
C_{\varepsilon_4}({\Theta})\,T \biggl( \delta ^{- (p + N(p - 1))}
+ \delta ^{- \frac{p + N(p - 1)}{p + N (p - 3)}}\biggr),
\end{multline}
$\forall\,\varepsilon_4 > 0,\ s \geqslant 0,\ \delta > 0,$ where
$C_{\varepsilon_4}({\Theta})$ is a constant depending on $p,N,
\varepsilon_4$ and $\Theta$. Choosing $\varepsilon_i
> 0,\ i =3,4$ sufficiently small, by the standard iteration procedure we
establish that
\begin{equation}\label{enql}
L_T(s_0 + \delta_0) \leqslant C({\Theta})\,T\,U(\delta_0)\
\forall\,s_0 \geqslant 0,\ \delta_0 > 0,
\end{equation}
where
$$
U(\delta_0):=  \delta_0^{- (p + N(p - 1))} +  \delta_0^{- \frac{p
+ N(p - 1)}{p + N (p - 3)}}.
$$
Let $\delta \to + \infty$ in (\ref{es1}), (\ref{es2}) and use
(\ref{enql}) for $s_0 = 0$ and $\delta_0 = s > 0$. We eventually
obtain
$$
A_T(s) \leqslant C({\Theta})T U^{\theta_1}(s),\ B_T(s) \leqslant
C({\Theta})T U^{3\theta_2/p}(s),
$$
for every $s > 0$. Hence (\ref{esde}) follows.

\end{proof}

\begin{proof}[Proof of Theorem \ref{Th2}]
From (\ref{J1}) and the decay estimations (\ref{esde}), we get
\begin{multline*}
G(s_0) := s_0 + \frac{1}{1 - \varepsilon}J_T(s_0) \leqslant s_0 +
c\, \max \{ F^{\frac{1}{p\beta}}(T)(A_T^{\frac{\beta_1}{p(1 +
\beta_1)}} (s_0)+ B_T^{\frac{\beta_1}{p(1 + \beta_2)}} (s_0)),\\
F^{\frac{1}{\beta}}(T)(A_T^{\frac{\beta_2}{1 + \beta_1}} (s_0) +
B_T^{\frac{\beta_2}{1 + \beta_2}} (s_0))\} \leqslant s_0 +
c(\Theta)\, \max \{ F^{\frac{1}{p\beta}}(T)T^{\frac{\beta_1}{p(1 +
\beta_1)}} s_0^{- \frac{N (p - 1)\beta_1}{p(1 + \beta_1)}},\\
F^{\frac{1}{p\beta}}(T) T^{\frac{\beta_1}{p(1 + \beta_2)}} s_0^{-
\frac{2N \beta_1}{p(p+N(p-3))(1 + \beta_2)}},\
F^{\frac{1}{\beta}}(T) T^{\frac{\beta_2}{1 + \beta_1}}
s_0^{-\frac{N \beta_2(p - 1)}{1 + \beta_1}},\\
F^{\frac{1}{\beta}}(T)T^{\frac{\beta_2}{1 + \beta_2}} s_0^{-
\frac{2N\beta_2}{(p+N(p-3))(1 + \beta_2)}}\} \leqslant
\tilde{G}(s_0) := s_0 + \\ c(\Theta)\, \max \{ T^{\frac{\beta_1 +
\alpha_1}{p(1 + \beta_1)}} s_0^{- \frac{N \beta_1(p - 1)}{p(1 +
\beta_1)}},  T^{\frac{\beta_2+ \alpha_2}{1 + \beta_2}} s_0^{-
\frac{2N\beta_2}{(p+N(p-3))(1 + \beta_2)}} \}.
\end{multline*}
Minimizing the function $\tilde{G}(s_{0}) $ we suppose that
\begin{multline}\label{rtr}
s = \Gamma (T) = \tilde{G}(s_{0, min}) = c({\Theta})\, \max \{
T^{\frac{\beta_1 + \alpha_1}{p(1 + \beta_1) + N \beta_1(p-1)}},\
T^{\frac{(\beta_2 + \alpha_2)(p+N(p-3))}{(1 + \beta_2)(p+N(p-3)) +
2N\beta_2}}\}= \\
c({\Theta})\, \max \{  T^{\frac{1}{p + N(p - 2)}},\ T^{\frac{p +
N(p - 3)}{p + N(p - 2)}}\}.
\end{multline}
 Thus $u_k(T,x) \equiv 0$ for all $x \in \{x=(x',x_N) : x_N
\geqslant \Gamma(T)\},\ k=\overline{1,N}$. And Theorem~1.2 is
proved completely.

\end{proof}

\section*{Appendix A}
\renewcommand{\thesection}{A}\setcounter{equation}{0}
\setcounter{lemma}{0}
\renewcommand{\thetheorem}{A.\arabic{lemma}}

\begin{lemma}$\cite{SH2}$
Let the nonnegative continuous nonincreasing function $f(s):
[s_0,\infty) \to \mathbb{R}^1$ satisfies the following functional
relation:
$$
f(s + f(s)) \leqslant \varepsilon\,f(s)\ \forall\, s \geqslant
s_0,\ 0 < \varepsilon < 1.
$$
Then $f(s)\equiv 0 \ \forall\, s \geqslant s_0+(1 -
\varepsilon)^{-1}f(s_0)$.
\end{lemma}

\begin{lemma}$\cite{Ni}$
If $\Omega  \subset \mathbb{R}^N $ is a bounded domain with
piecewise-smooth boundary, $a > 1$, $b \in (0, a),$ and $d
> 1,$ then there exist positive constants
$d_1$ and $d_2$ $(d_2 = 0$ if the domain $\Omega$ is unbounded and
$d_2 = \tilde{d}_2 \delta^{ - \frac{N(a-b)}{ab}}$ if the domain
$\Omega = K(s,\delta))$ that depend only on $\Omega ,\ d,\ b,$ and
$N$ and are such that, for any function $v(x) \in W^{1,d} (\Omega
) \cap L^b (\Omega )$, the following inequality is true:
$$
\left\| { v} \right\|_{L^a (\Omega )}  \leqslant d_1 \left\| {D v}
\right\|_{L^d (\Omega )}^\theta  \left\| v \right\|_{L^b (\Omega
)}^{1 - \theta }  + d_2 \left\| v \right\|_{L^b (\Omega )}
$$
where $\theta  = \frac{{\tfrac{1} {b}  - \tfrac{1} {a}}}
{{\tfrac{1} {b} + \tfrac{1} {N} - \tfrac{1} {d}}} \in \left[ {0,1}
\right)$.
\end{lemma}

{\footnotesize {\bf Acknowledgement.}We would like to thank to
Prof. \v Sarka Ne\v casov\' a and Prof. Eduard Feireisl for they
valuable comments and remarks.}

\smallskip
$$
\begin{gathered}
\text{\normalsize\it Yu.V. Namlyeyeva}\hfill\\
\text{\normalsize\it Institute of Applied Mathematics}\hfill\\
\text{\normalsize\it and Mechanics of NAS of Ukraine},\hfill\\
\text{\normalsize\it R. Luxemburg str. 74, 83114 Donetsk,}\hfill\\
\text{\normalsize\it Ukraine}\hfill\\
\text{\normalsize\it e-mail: namleeva@iamm.ac.donetsk.ua}\hfill\\
\end{gathered}
$$
\smallskip
$$
\begin{gathered}
\text{\normalsize\it R.M. Taranets}\hfill\\
\text{\normalsize\it Institute of Applied Mathematics}\hfill\\
\text{\normalsize\it and Mechanics of NAS of Ukraine},\hfill\\
\text{\normalsize\it R. Luxemburg str. 74, 83114 Donetsk,}\hfill\\
\text{\normalsize\it Ukraine}\hfill\\
\text{\normalsize\it e-mail: taranets\_r@iamm.ac.donetsk.ua}\hfill
\end{gathered}
$$


\begin{thebibliography}{99}

\bibitem{A} Antontsev, S.N., Diaz, J.I., Shmarev, S. Energy methods for free boundary
problems: applications for nonlinear PDEs and fluid mechanics
(Progress in Nonlinear Differential Equations and Their
Applications) \textbf{XI}, A Birkh\"{a}user book, 2002, 329 p.

\bibitem{AO1} Antontsev, S.N., de Oliveira, H.B. Navier-Stokes
equations woth absorbtion under slip boundary conditions:
existence, uniqueness and extinction in time \emph{RIMS
K\^{o}ky\^{u}roku Bessantsu} \textbf{B1} (2007), Kyoto University,
21--42

\bibitem{AO2} Antontsev, S.N., de Oliveira, H.B. Localization of
weak solutions for non-Newtonian fluid flows (Portuguese)
\emph{Proceedings of the Congress of Computational Methods in
Engineering}, National Labaratory of Civil Engeneering, Lisbon,
(2004), 15 p.

\bibitem{AT} Andreucci, D., Tedeev, A. A Fujita type result for a degenerate
Neumann problem in domains with noncompact boundary. \emph{J.
Math. Anal. Appl.} \textbf{231} (1999), 543--567

\bibitem{BBN}
Bellout, H., Bloom, F., Ne\v cas, J.: Solutions for
incompresssible non-Newtonian fluids, \emph{C. R. Acad. Sci.
Paris} \textbf {317}, S\'erie I, (1993), 795-800

\bibitem{KV} Kamin, S., Vazquez, J. L., Fundamental solutions and
asymptotic behavior for the p-Laplacian equation, \emph{Rev. Mat.
Iberoamericana} \textbf{4}(2), (1988), 339--354

\bibitem{LPV} Lee, K., Petrosyan, A., Vazquez, J. L., Large time
geometric properties of solutions of the evolution p-Laplacian
equation, \emph{Journal of Diff. Equations} \textbf{229}, (2006),
389--411

\bibitem{MNRR}
M\' alek, J., Ne\v cas, J., Rokyta, M., R\accent23 u\v zi\v cka,
M., Weak and Measure-valued Solutions to Evolutionary PDEs,
\emph{Applied Mathematics and Mathematical Computation}
\textbf{13}, Chapman $\&$ Hall, 1994

\bibitem{Ni} Nirenberg, L. , An extended interpolation
inequality, \emph{Ann. Scuola Norm. Sup. Pisa} \textbf{20} (1966),
733--737

\bibitem{Po1}
Pokorn\' y, M., Cauchy problem for the non-Newtonian
incompressible fluid,\emph{Applications of Mathematics},
\textbf{41}, No.3, (1996), 169-201

\bibitem{Po2}
Pokorn\' y, M., Cauchy problem for the non-Newtonian
incompressible fluid (Master thesis), Faculty of Mathematics and
Physics, Charles University, Prague, 1996

\bibitem{Ra}
Rajagopal, K.R., \emph{Mechanics of non-Newtonian fluids}, Ed. G.
P. Galdi, J. Ne\v cas : \emph{Recent Developments in Theoretical
Fluid Dynamics. Pitman Research Notes in Math.series}
\textbf{291}, Longman Scientific $\&$ Technical, Essex, 1993,
129-162

\bibitem{Rav} Raviart, P.A., Sur la r\' esolution de certaines \'
equations paraboliques non lin\' eaires, \emph{J. Functional
Analysis} \textbf{5} (1970), 299--328

\bibitem{SS}  Sapronov, D., Shishkov, A. Asymptotic behaviour
of supports of solutions of quasilinear many-dimensional parabolic
equations of non-stationary diffusion-convection type. \emph{Sb.
Math.} \textbf{197} (2006), 753--790

\bibitem{SSH} Shishkov, A. Dynamics of the geometry of the support of the
generalized solution of a higher-order quasilinear parabolic
equation in divergence form. \emph{Differential Equations}
\textbf{29} (1993), 460--469

\bibitem{SH2}  Shishkov,  A.,  Shchelkov, A. Dynamics of the
supports of energy solutions of mixed problems for quasi-linear
parabolic equations of arbitrary order. \emph{Izvestiya RAN: Ser.
Math.} \textbf{62} (1998), 601--626

\bibitem{T5}  Taranets, R. Propagation of perturbations in
thin capillary film equations with nonlinear diffusion and
convection. \emph{Siberian Math. J.} \textbf{47} (2006), 914--931

\bibitem{T3} Taranets, R.,  Shishkov,  A. Effect of time delay of
support propagation in equations of thin films. \emph{Ukrainian
Math. J.} \textbf{55} (2003), 1131--1152

\bibitem{Va1} Vazquez, J.L. The porous medium equation: mathematical
theory, Oxford: Clarendon, 2007

\bibitem{Va2} Vazquez, J.L. Smoothing and Decay estimates for
Nonlinear diffusion equationss, Equations of Porous medium type,
Oxford Lecture series in Mathematics and its applications
\textbf{33}, Oxford Press, 2006

\end{thebibliography}
\end{document}